\documentclass[12pt]{article}
\usepackage{floatflt}
\usepackage[english]{babel}
\usepackage{amsmath, amssymb}  
\usepackage{mathrsfs}
\usepackage{authblk}

\newtheorem{theorem}{Theorem}

\def\Pb{\mathbf{P}}
\def\Ex{\mathbf{E}}

\renewcommand{\P}{\mathbb P}

\newcommand{\II}{\mathbb I}
\newcommand{\BB}{\mathbb B}
\newcommand{\E}{\mathbb E}
\newcommand{\R}{{\cal R}}

\newcommand{\VV}{\mathbb V}

\def\1{\mbox{1\hspace{-.25em}I}}
\let\scr\mathscr

\begin{document}
\title{On APF  Test 
for Poisson Process with Shift and Scale Parameters}
\date{}

\date{}
\author[1]{A. S. Dabye}
\author[2]{Yu. A. Kutoyants}
\author[3]{E.D. Tanguep}
\affil[1,3]{\small Universit\'e Gaston Berger,  Saint--Louis, S\'en\'egal}
\affil[2]{\small  Le Mans University, Le Mans, France}
\affil[2]{\small   National Research University ``MPEI'', Moscow, Russia}
\affil[2]{\small  Tomsk State University, Tomsk, Russia}
\maketitle

\begin{abstract}
We propose the goodness of fit test for  inhomogeneous Poisson processes
with unknown scale and shift parameters. A test statistic of Cram\'er-von Mises type
is proposed and its asymptotic behavior is  studied.  We show that under null 
hypothesis the limit distribution of this statistic does not depend on unknown
parameters. 
\end{abstract}

\bigskip{}
\bigskip{}


\section{Introduction}

The problems of the construction of goodness of fit tests in the case of
i.i.d. observations are well studied \cite{LR05}. Special attention is payed
to the case of parametric null hypothesis. Wide class of distributions can be
parametrized by the shift and scale parameters, say, $F\left(\frac{x-\vartheta
  _1}{\vartheta _2}\right)$. In the case of such families several authors
showed that the limit distributions of the Kolmogorov-Smirnov and Cramer-von
Mises tests statistics do not depend on the unknown parameters (see
\cite{Da58}, \cite{Du73}, \cite{G53}, \cite{DN73}, \cite{Mart79}, \cite{Mart10} and
references therein). We call such tests {\it asymptotically parameter free} (APF).

For the continuous time stochastic processes the goodness of fit testing is
not yet well developed. We can mention here several works for diffusion and
Posson processes \cite{D13}, \cite{DTT16}, \cite{DK07}, \cite{Da77},
\cite{KK14}, \cite{Kut14},\cite{Kut14a}, \cite{W75}.  The problem of goodness
of fit testing for inhomogeneous Poisson process is interesting because there
is a wide literature on the applications of inhomogeneous Poisson process
models in different domains (astronomy, biology, image analysis, medicine,
optical communication, physics, reliability theory, etc.).  Therefore to know
if the observed Poisson process corresponds to some parametric family of
intensity functions is  important. 

We consider the problem of goodness of fit testing for inhomogeneous Poisson
process which  under the null hypothesis  has the intensity function with
shift and scale parameters. We show that as in the classical case the limit
distribution of the Cramer-von Mises type statistics does not depend  on these
unknown parameters. This allows us to construct the corresponding APF goodness of
fit test of fixed asymptotic size. 

\section{Statement of the problem and auxiliary results}

Suppose that we observe $n$ independents inhomogeneous Poisson processes
$X^{n}=\big(X_1,\ldots,X_n\big)$, where $X_j=\Bigl(X_j\left(t\right), t\in
\R\Bigr)$ are trajectories of the Poisson processes with the mean function 
$
\Lambda
\left(t\right)=\E X_j\left(t\right)=\int_{-\infty }^{t}\lambda
\left(s\right){\rm d}s. 
$
Here $\lambda \left(\cdot \right)\geq 0$ is the corresponding intensity
function. 

 Let us remind the construction of GoF test of Cram\'er-von Mises type in the
 case of simple null hypothesis. The class of tests $\left(\bar\Psi
 _n\right)_{n\geq 1}$ of asymptotic size $\varepsilon \in (0,1)$ is
\begin{align*}
{\scr K}_\varepsilon =\left\{\bar\Psi _n\;:\qquad \lim_{n\rightarrow \infty
}\E_0\bar\Psi _n=\varepsilon  \right\}.
\end{align*}
 
Suppose  that  the basic hypothesis
is simple, say,
$
{\scr H}_0 ~:~~ \Lambda \left(\cdot \right) = \Lambda_0 \left(\cdot\right),
$
where $\Lambda_0 \left(\cdot  \right)$ is a know continuous function
satisfying 
$\Lambda_0\left(\infty \right) < \infty $. The alternative is composite (non
parametric) 
$
{\scr H}_1 ~:~~ \Lambda \left(\cdot \right)\not =\Lambda_0 \left(\cdot
\right).
$
Then we can introduce the Cram\'er-von Mises (C-vM) type statistic 
$$
\tilde \Delta _n = \frac{n}{\Lambda_0 \left(\infty
  \right)^2}\,\int_{\R}^{}\left[\widehat{\Lambda}_{n}(t) -
 \Lambda_0(t)\right]^{2}\,{\rm d}\Lambda_0(t),
$$
where 
$
\widehat{\Lambda}_{n}(t)=\frac{1}{n}\sum_{j=1}^{n}X_j\left(t\right)
$ is the empirical mean of the Poisson process. It can be verified that under
${\scr H}_0$ this statistic converges to the following limit:
$$
\tilde \Delta _n \Longrightarrow \Delta \equiv  
\int_{0}^{1}W\left(s\right)^2{\rm d}s, 
$$
where $W\left(s\right),\; 0 \leq s \leq 1 $  is a standard Wiener process.
Therefore the C-vM type test $\tilde\psi_n\left(X^n\right) =
1\!\!1_{\left\{\tilde \Delta _n > c_\varepsilon \right\}}$ with the threshold
$c_\varepsilon $ defined by the equation $ \mathbb{P} \left\{\Delta >
c_\varepsilon \right\} = \varepsilon$ belongs to ${\scr K}_\varepsilon $.
 This test is {\it asymptotically distribution free} (ADF) (see, e.g.,
 \cite{DK07}). Remind that   the test is called ADF  if the limit 
distribution of the test statistic under hypothesis does not depend on the
mean function  $\Lambda _0\left(\cdot \right)$. 

Let us consider the case of the parametric null hypothesis. It can be
formulated as follows. We have to test the null hypothesis 
$$
 {\scr H}_0 ~:~~ \Lambda \left(\cdot \right)\in {\cal L }\left(\Theta \right) = 
\Bigl\{\Lambda_0 \left(\vartheta,t\right), ~ \vartheta \in \Theta, ~ t\in
\R \Bigr\},  
$$
against the alternative
$
{\scr H}_1 ~: \Lambda \left(\cdot \right) \notin {\cal L }\left(\Theta \right).
$
Here $\Lambda_0(\vartheta,\cdot)$ is a known mean function of the Poisson
process depending on some finite-dimensional unknown parameter $\vartheta \in
\Theta \subset {\R}^d$. Note that  under ${\scr H}_0 $ there exists the {\it
  true value} $\vartheta _0\in\Theta $ such that the mean of the observed
Poisson process $\Lambda \left(t\right)=\Lambda \left(\vartheta _0,t\right),
t\in \R$.

The C-vM  type GoF  test can be constructed by a similar way. Introduce the
normalized process 
 $ \bar{ u}_n(t)\equiv u_n\left(t, \bar{\vartheta}_n \right)= \sqrt{n}\left(
  \widehat{\Lambda}_{n}(t)- 
\Lambda_0( \bar{\vartheta}_n,t) \right),\; t\in  \R.$
Here $\bar\vartheta _n$ is some estimator of the parameter $\vartheta $, which
is (under hypothesis ${\scr H}_0$) consistent and asymptotically normal $\sqrt{n}\left(\bar\vartheta _n-\vartheta _0 \right)\Longrightarrow
\xi $.

The corresponding C-vM type statistic can be 
\begin{align*}
\bar\Delta _n=\frac{n}{ \Lambda _0\left(\infty
  ,\bar{\vartheta}_n\right)^2}\int_{\R }^{}\left(
  \widehat{\Lambda}_{n}(t)- 
\Lambda_0(\bar{\vartheta}_n,t) \right)^2{\rm d} \Lambda_0( \bar{\vartheta}_n,t) 
\end{align*}

Then, under null hypothesis ${\scr H}_0$, we can verify  the convergence
\begin{align*}
\bar{ u}_n(t)&=\sqrt{n}\left(
  \widehat{\Lambda}_{n}(t)- \Lambda_0(
 \vartheta_0,t)\right)+\sqrt{n}\left(\Lambda_0(\vartheta_0,t ) -
\Lambda_0(\bar{\vartheta}_n,t) \right)\\
&=W_n\left(t
\right)-\langle\sqrt{n}\left(\bar{\vartheta}_n-\vartheta _0 \right),\frac{\partial
  \Lambda_0(\vartheta_0,t )}{\partial \vartheta }\rangle +o\left(1\right) \\
&\Longrightarrow W\left(\Lambda_0(\vartheta_0,t)\right)-\langle\xi\left(\vartheta
_0\right),\dot \Lambda _0\left(\vartheta_0,t\right)\rangle.
\end{align*}
Here $ \langle,\rangle$ is the scalar product in ${\R}^d$ and dot means
differentiation w.r.t. $\vartheta $. Let us denote $s=\Lambda _0\left(\vartheta _0,\infty
\right)^{-1}\Lambda _0\left(\vartheta _0,t\right)$ and introduce the vector 
$
G\left(\vartheta _0,s\right) =\Lambda _0\left(\vartheta _0,\infty
\right)^{-1/2}\dot \Lambda _0\left(\vartheta_0,t\right) .
 $
Then we obtain the convergence
\begin{align*}
\bar\Delta _n\Longrightarrow \bar\Delta\left(\vartheta _0,\Lambda
_0\right)=\int_{0}^{1}\left[W\left(s\right)-\langle\xi\left(\vartheta 
_0\right),G\left(\vartheta _0,s\right)\rangle \right]^2{\rm d}s,
\end{align*}
where $W\left(s\right),0\leq s\leq 1$ is standard Wiener process.  Here the
distribution of the limit random variable $\bar\Delta\left(\vartheta
_0,\Lambda _0\right) $ depends on the true value $\vartheta _0$ and on the
mean function $\Lambda _0\left(\vartheta _0,\cdot \right)$.

Therefore if we propose a GoF test based on this statistics, say, $\Phi
_n=\1_{\left\{\bar \Delta _n>c_\varepsilon \right\}} $,  then to find the
threshold $c_\varepsilon $ such that $\Phi
_n\in {\scr K}_\varepsilon   $ we have to solve the equation
$\P_{\vartheta _0}\left(\bar\Delta\left(\vartheta _0,\Lambda _0\right)
>c_\varepsilon \right)=\varepsilon$. 
 The solution $c_\varepsilon =c_\varepsilon \left(\vartheta _0,\Lambda
 _0\right)$, where $\vartheta _0$ is the unknown true value. There are several
 possibilities to construct the test belonging ${\cal K}_\varepsilon $. One is
 to calculate the function $c_\varepsilon \left(\vartheta _0,\Lambda
 _0\right) $, verify that this function is continuous w.r.t. $\vartheta $ and
 then to use the consistent estimator $\bar\vartheta _n$  for the threshold
 $c_\varepsilon \left(\bar\vartheta _n,\Lambda  _0\right)  $. Another
 possibility is to use the linear transformation of the statistic $\bar
 u_n\left(\cdot \right)$, which transforms it in the Wiener process (see, e.g.,
 \cite{Kh81} or \cite{KK14}). In this work we follow the third approach: we show
 that the limit distribution of the statistic  does not
 depend on $\vartheta _0$. 

In particular, the goal of this work is to show that if the unknown parameter is
two-dimensional $\vartheta =\left(\alpha ,\beta \right)$, where $\alpha \in \R$ is
the {\it shift} and $\beta \in\R_+$ is the {\it scale} parameters, then it is possible to
construct a test statistic  $\hat\Delta _n$   whose limit distribution does not depend on
$\vartheta _0$.  The mean function under null hypothesis is 
\begin{align*}
\Lambda _0\left(\vartheta,t \right)=\int_{-\infty }^{t}\lambda
_0\left(\frac{v-\alpha }{\beta }\right){\rm d}v,\qquad t\in\R. 
\end{align*}
The proposed test statistic is
\begin{align*}
\hat\Delta _n=\frac{n}{\hat\beta _n^2}\int_{\R}^{}\left[\hat\Lambda _n\left(t\right)-\Lambda
  _0(\hat\vartheta _n,t)\right]^2{\rm d}\Lambda
_0(\hat\vartheta _n,t) .
\end{align*}
Here $\hat\vartheta _n$ is the maximum likelihood estimator (MLE) of the
vector parameter $\vartheta $. We show that $\hat\Delta _n\Rightarrow \Delta
$, where $\Delta =\Delta \left(\Lambda _0\right)$, i.e., the distribution of
the random variable $\Delta\left(\Lambda _0\right) $ does not depend on
$\vartheta _0$. Remind that the function $\Lambda _0\left(t \right),t\in\R $ is
known and therefore the solution $c_\varepsilon =c_\varepsilon \left(\Lambda
_0\right)$ can be calculated before the experiment using, say, numerical
simulations.

We are given $n$ independent observations $X^n=\left(X_1,\ldots,X_n\right)$ of
inhomogeneous Poisson processes $X_j=\left(X_j\left(t\right),t\in\R \right)$
with the mean function $\Lambda \left(t\right)=\E X_j\left(t\right),t\in\R$.
 We have to construct a GoF test in the
hypothesis testing problem with parametric null hypothesis ${\scr H}_0$. More
precizely, we suppose that under  ${\scr H}_0$ the mean function $\Lambda
\left(t \right)$ is absolutely continuous: $\Lambda '\left(t\right)= \lambda
_0\left(\vartheta _0,t\right)$. 
Here $\vartheta _0$ is the true value and  the intensity function is
$\lambda _0\left(\vartheta_0 ,t\right)=\lambda _0\left(\frac{t-\alpha_0
}{\beta_0 }\right),\;\vartheta 
_0=\left(\alpha _0,\beta _0\right)\in\Theta \subset {\R}^2. 
$
The set $\Theta= \left(a_1,a_2\right)\times (b_1,b_2)$ and $b_1>0$, where all
constants are finite.
Therefore if we denote 
$
\Lambda_0
\left(t\right)=\int_{-\infty }^{t}\lambda _0\left(v\right){\rm d}v,\; t\in\R,
$
 then the mean function under null hypothesis is 
\begin{align*}
\Lambda \left(t\right)=\Lambda_0 \left(\vartheta _0,t\right)=\beta_0 \Lambda _0\left(\frac{t-\alpha _0}{ \beta _0}\right).
\end{align*}
It is convenient to use two different functions $\Lambda _0\left(\vartheta,t
\right)$ and $\Lambda _0\left(t\right)$ and we hope that such notation will
not be misleading. 

Therefore, we  have the parametric null hypothesis
\begin{align*}
{\scr H}_0:\qquad \quad \Lambda \left(\cdot \right)\in {\cal L }\left(\Theta \right),
\end{align*}
where the  parametric family is 
\begin{align}
\label{h0}
&{\cal L }\left(\Theta \right)=\Bigl\{\Lambda \left(\cdot \right):\Lambda \left(t\right) =
  \beta\Lambda_0\left(\frac{t - \alpha}{\beta}\right),\;t\in{\R},\quad  \vartheta =
  (\alpha,\beta) \in \Theta \Bigr\}.
\end{align}
Here $\Lambda_0 \left(\cdot \right)$ is a known absolutely continuous function
with properties:\\ $
  \Lambda_0 \left(-\infty\right)=0,\;\Lambda_0 \left(\infty\right)<\infty 
$.

 We consider the class of tests of asymptotic level $ \varepsilon $:
\begin{align}
\label{kepsil}
\mathcal{K}_{\varepsilon}=\left\{\bar\Psi_n: \lim_{n\rightarrow \infty
}\E_\vartheta\bar\Psi_n  =\varepsilon,\quad \vartheta \in \Theta  \right\}.
\end{align}
The test studied in this work is based on the following  statistic of C-vM type:
\begin{align}
\label{kepsil1}
\hat\Delta _n =
\frac{n}{\hat{\beta}_n^2}\,\int_{\R}^{}\left[\hat{\Lambda}_{n}(t)-\hat\beta _n
  \Lambda_0   \left(\frac{t-\hat{\alpha }_n}{\hat\beta _n} \right)
  \right]^{2}\,\lambda_0\left(\frac{t-\hat{\alpha }_n}{\hat\beta _n} \right)
\;{\rm d}t .
\end{align}
where $\hat{\vartheta}_n=(\hat\alpha _n,\hat\beta _n)$ is the MLE. Remind that 
the log-likelihood ratio for this model of observations is
\begin{align*}
\ln L\left(\vartheta ,\vartheta_1 ,X^n\right)=\sum_{j=1}^{n}\int_{\R}^{}\ln
\frac{\lambda_0 \left(\vartheta,t \right)}{\lambda_0 \left(\vartheta,t_1 \right)} {\rm
  d}X_j\left(t\right)-n\int_{\R}^{}\left[\lambda_0 \left(\vartheta,t
  \right)-\lambda_0 \left(\vartheta,t_1 \right)\right]{\rm d}t ,
\end{align*}
and the MLE $\hat\vartheta _n$ is defined by the equation
\begin{align}
\label{mle}
L\left(\hat\vartheta _n,\vartheta_1 ,X^n\right)=\sup_{\vartheta \in\Theta
}L\left(\vartheta,\vartheta_1 ,X^n\right). 
\end{align}
 Here $\vartheta _1\in\Theta $ is some fixed value.

As we use the asymptotic properties of the MLE $\hat\vartheta_n$, we need some
regularity conditions, which we borrow from \cite{Kut98} (see the conditions
{\bf B1}-{\bf B5} in the Section 2.1 there). 

Note that the derivative (vector) $\frac{\partial \lambda \left(\vartheta
  ,t\right)}{\partial \vartheta } $   of the intensity function is 
\begin{align}
\label{def}
\dot\lambda \left(\vartheta ,t\right)=\left(\frac{\partial \lambda
  \left(\vartheta ,t\right)}{\partial \alpha }, \frac{\partial \lambda
  \left(\vartheta ,t\right)}{\partial \beta  }   \right)=-\lambda '\left(
\frac{t-\alpha }{\beta }\right)  \;\left(\frac{1}{\beta }\;,\frac{t-\alpha
}{\beta ^2}\right) .
\end{align}
Here $\lambda '\left(t\right)=\frac{{\rm d}\lambda \left(t\right)}{{\rm
    d}t}$.

 {\it Conditions }${\scr R}$

${\scr R}_1$. {\it The intensity function $\lambda_0\left(\cdot\right)$ is
   strictly positive and two times continuously differentiable.}

${\scr R}_2$. {\it For any $\vartheta _0\in\Theta $ we have}
\begin{align}
\label{c1}
&\lim_{\left\|\vartheta -\vartheta _0\right\|\rightarrow 0}\int_{\R}^{}\left|
\frac{\dot \lambda_0 \left(\vartheta ,t\right)}{\sqrt{\lambda_0 \left(\vartheta
    ,t\right)}} -\frac{\dot \lambda_0 \left(\vartheta_0 ,t\right)}{\sqrt{\lambda_0 \left(\vartheta_0
    ,t\right)}}\right|^2\lambda_0 \left(\vartheta _0,t\right)\,{\rm d}t=0,\\
&\qquad \qquad \sup_{\vartheta \in\Theta }\int_{\R}^{}\left|
\frac{\dot \lambda_0 \left(\vartheta ,t\right)}{\sqrt{\lambda_0 \left(\vartheta
    ,t\right)}} \right|^4\lambda_0 \left(\vartheta _0,t\right)\,{\rm
  d}t<\infty \label{c2} .
\end{align}

${\scr R}_3$. {\it The function $\lambda _0\left(\cdot \right)$ satisfies the conditions}
\begin{align}
\label{c3}
\int_{\R}^{}t^2\lambda _0\left(t\right)\;{\rm d}t<\infty ,\qquad \int_{\R}^{}t^4\left|\lambda _0'\left(t\right)\right|\;{\rm d}t<\infty .
\end{align}

Of course, we suppose that the expressions under the sign of integrals are
integrable in the required sense. 

For the consistency of the MLE we need the identifiability condition

${\it I}.$ {\it For any $\nu >0$}
\begin{align*}
\inf_{\left\|\vartheta -\vartheta _0\right\|>\nu }
\int_{\R}^{}\left[\sqrt{\lambda_0 \left(\vartheta ,t\right)}-\sqrt{\lambda_0
    \left(\vartheta_0 ,t\right)} \right]^2{\rm d}t>0.
\end{align*}

Note that in the case of shift and scale parameters this condition is
fulfilled. Indeed, suppose that for some $\nu >0$ this integral is 0. Then there exists
$\vartheta _1\not=\vartheta _0$ ($\left\|\vartheta _1-\vartheta _0\right\|\geq
\nu $) such that $\lambda \left(\vartheta_1 ,t\right)\equiv \lambda
\left(\vartheta_0 ,t\right) $. Recall that the functions are
continuous. Therefore 
$\lambda \left(\frac{t-\alpha _1}{\beta _1}\right)=\lambda \left(\frac{t-\alpha
  _0}{\beta _0}\right),\;  t\in \R$
or after the change of variables $s=\beta _0^{-1}\left(t-\alpha _o\right)$ we have
\begin{align*}
\lambda _0\left(s\right)=\lambda _0\left( \frac{\beta _0}{\beta _1}\;
s-\frac{\alpha _1-\alpha _0}{\beta _1}\right),\qquad  s\in \R. 
\end{align*}
Of course, such function  $\lambda _0\left(\cdot \right)\not \in L_1\left(\cal
R\right)$.  Hence, the condition of identifiability is fulfilled. 

To construct the test statistics we need the following property of the mean
function

 {\it For all $\vartheta _0\in\Theta $}
\begin{align}
\label{c4}
\sup_{\vartheta \in\Theta }\int_{\R}^{}\left|  \dot \Lambda_0 \left(\vartheta ,t\right)\right|^2\lambda
\left(\vartheta _0,t\right)\,{\rm d}t <\infty .
\end{align}
This condition   can be expressed in terms of the
 function  $\lambda _0\left(\cdot \right)$ like \eqref{c1}-\eqref{c2}. Indeed
 we have
\begin{align*}
\left|  \dot \Lambda_0 \left(\vartheta ,t\right)\right|^2=\lambda _0\left( \frac{t-\alpha }{\beta}\right)
^2+\left|\Lambda _0\left( \frac{t-\alpha
}{\beta}\right)-\left(\frac{t-\alpha }{\beta }\right)\lambda _0\left(
\frac{t-\alpha }{\beta}\right) \right|^2 .
\end{align*}
As the function $\lambda _0\left(\cdot \right)$ is bounded, it is sufficient
to suppose \eqref{c3} and we obtain \eqref{c4}.

Let us introduce the Fisher information matrix
\begin{align*}
\II\left(\vartheta \right)=\frac{1}{\beta}
\left(\begin{array}{c c c}
\int_{\R}
\frac{\lambda'_0 (s)^{2}}{\lambda_0(s)} \,{\rm d}s&  &
\int_{\R}\frac{s \;\lambda'_0(s)^{2}}{\lambda_0(s)} \,{\rm d}s\\
             &  &                 \\
\int_{\R}\frac{s \,\lambda'_0(s)^{2}}{\lambda_0(s)} \,{\rm d}s &  &
\int_{\R}\frac{s^{2}\,\lambda'_0(s)^{2}}{\lambda_0(s)} \, {\rm d} s\\
\end{array}\right)=\frac{1}{\beta } \;\II_*,
\end{align*}
where the matrix $\II_* $ does not depend on $\vartheta $.
Note that the matrix  $\II_* $ is non degenerate. Indeed, the determinant is
\begin{align*}
D=\int_{\R}
\frac{\lambda'_0 (s)^{2}}{\lambda_0(s)} \,{\rm
  d}s\int_{\R}\frac{s^{2}\,\lambda'_0(s)^{2}}{\lambda_0(s)} \, {\rm d} s
-\left( \int_{\R}\frac{s \;\lambda'_0(s)^{2}}{\lambda_0(s)} \,{\rm d}s  \right)^2.
\end{align*}
Remind that by Cauchy-Schwartz inequality
\begin{align*}
\left( \int_{\R}\frac{s \;\lambda'_0(s)^{2}}{\lambda_0(s)} \,{\rm d}s
\right)^2\leq \int_{\R} 
\frac{\lambda'_0 (s)^{2}}{\lambda_0(s)} \,{\rm
  d}s\int_{\R}\frac{s^{2}\,\lambda'_0(s)^{2}}{\lambda_0(s)} \, {\rm d} s.
\end{align*}
The equality in  Cauchy-Schwartz inequality ($D=0$) we obtain if and only if
$\left|s\lambda _0'\left(s\right)\right|\equiv \left|\lambda _0'\left(s\right)\right|,\; s\in\R.$
Of course such equality is impossible, if $\lambda '\left(s\right)\not=0$ or
$s\not= \pm 1$. 
As the function $\lambda_0 \left(\cdot \right)$ is positive and differentiable,
 we have
\begin{align*}
\int_{\R}\frac{\lambda'_0(s)^{2}}{\lambda_0(s)} \,{\rm d}s >0. 
\end{align*}

We suppose that the intensity function  $\lambda _0\left(t\right),t\in \R$  is
strictly positive because if we have 
a set of positive Lebesgue measure, where $\lambda _0\left(t\right)=0$ and the
unknown parameters are shift and scale, then the measures induced by the
observations will be not equivalent. The properties of the MLE will be
different.

 Under these conditions, the MLE is uniformly 
consistent, asymptotically normal 
$\sqrt{\frac{n}{\beta _0}}\left(\hat\vartheta _n-\vartheta
_0\right)\Longrightarrow \zeta \sim {\cal N}\left(0, \II_*^{-1}\right)$
and the polynomial moments converge 
\begin{align}
\label{cm}
\lim_{n\rightarrow \infty } \left(\frac{n}{\beta _0}\right)^\frac{p}{2}\Ex_{\vartheta _0}\|\hat\vartheta
_n-\vartheta _0  \|^p =\Ex\left\|\zeta \right\|^p.
\end{align}
For the proof see Theorem 2.4 in \cite{Kut98}. Note that the distribution of
the vector $\zeta $ does not depend on $\vartheta_0 $.

\section{Main result}

 Introduce the following random variable:
\begin{align}
\label{stati}
\Delta_0 & = \int_{\R} \left[W\big(\Lambda_0 \left(t\right)\big) -\langle
  \zeta ,\dot \Lambda _0\left(t\right)\rangle\right]^2 \,{\rm d}\Lambda _0\left(t\right),
\end{align}
where  $\dot \Lambda \left(\vartheta ,t\right)= \Bigl(-\lambda _0\left(s\right),
\Lambda_0\left(s\right)-s\lambda _0\left(s\right)\Bigr) $ and  $W\left(r \right),0\leq r\leq \Lambda _0\left(\infty \right)$ is a Wiener process.
The main result of this work is the following theorem.

\begin{theorem}
\label{T1} Let the conditions ${\cal R}$ be fulfilled then the test 
\begin{align*}
\hat \Psi _n\left(X^n\right)=1\!\!1_{\big\{\Delta _n>c_\varepsilon
  \big\}},\qquad \Pb \left(\Delta_0>c_\varepsilon \right)=\varepsilon 
\end{align*}
belongs to the class ${\scr K}_\varepsilon  $.
\end{theorem}
{\bf Proof.} We can write  
\begin{align*}
 \hat u_n\left(t\right)&=\sqrt{n} \left(\hat\Lambda _n\left(t\right)-\hat\beta _n
  \Lambda_0   \left(\frac{t-\hat{\alpha }_n}{\hat\beta _n} \right) \right)\\
&=\sqrt{n} \left(\hat\Lambda _n\left(t\right)- \Lambda _0\left(\vartheta_0,t
  \right)  \right)+ \sqrt{n}  \left(\Lambda _0\left(\vartheta_0,t \right) - 
  \Lambda_0   (\hat\vartheta _n,t ) \right)\\
&=W_n\left(t\right)-\langle\sqrt{n}\left( \hat\vartheta _n-\vartheta
  _0\right),\dot \Lambda _0\left(\vartheta_0,t
  \right)\rangle+r_n\left(t\right)\equiv u_n\left(t\right)+r_n\left(t\right). 
\end{align*}
Here the vector $\dot \Lambda _0\left(\vartheta,t \right)=\frac{\partial
}{\partial\vartheta } \Lambda _0\left(\vartheta,t \right)$ and we used the Taylor
formula.

We have to show that under the null hypothesis 
\begin{align}
\label{15}
&\frac{1}{\hat\beta _n^2}\int_{\R}^{}u_n\left(t\right)^2\lambda _0(\hat\vartheta _n,t)\,{\rm d}t \Longrightarrow
\int_{\R}^{}\left[W\left(\Lambda _0\left( s\right)\right)+\langle \zeta
  ,V\left(s\right)\right]^2 {\rm d}\Lambda _0\left(s\right) ,\\
\label{16}
&\int_{\R}^{}r_n\left(t\right)^2\lambda _0(\hat\vartheta _n,t)\,{\rm d}t\longrightarrow 0.
\end{align}
Here  $V\left(s\right)=\Bigl(\lambda _0\left(s\right),s\lambda
_0\left(s\right)-\Lambda _0\left(s\right)\Bigr)  $.

The convergences \eqref{15}, \eqref{16} we will prove in several steps. 
\begin{description}
\item[A].   We show that we have the convergence of finite dimensional distributions
\begin{align}
\label{fdd}
\left(\hat\beta _n^{-1} \hat u_n\left(t_1\right),\ldots, \hat\beta _n^{-1} \hat
u_n\left(t_k\right) \right)\Longrightarrow \left(\hat u\left(s_1\right),
\ldots,  \hat u\left(s_k\right)\right),
\end{align}
where we put $s_i=\beta _0^{-1}\left(t_i-\alpha _0\right)$ and 
$\hat u\left(s\right)=W\left(\Lambda _0\left( s\right)\right)+\langle \zeta
  ,V\left(s\right)\rangle .$

\item[B]. We verify the estimate: for $\left|t_1\right|<L,\left|t_2\right|<L$
  and any $L>0$
\begin{align}
\label{dd}
\Ex_{\vartheta _0}\left| \hat u_n\left(t_1\right)- \hat u_n\left(t_2\right)
\right|^2\leq C\left(1+L\right)\left|t_1-t_2\right|, 
\end{align}
where the constant $C>0$ does not depend on $n$.

\item[C].  We show that for any $\delta  >0$ there exists $L>0$ such that for all $n$
\begin{align*}
\int_{\left|t\right|>L}^{}\Ex_{\vartheta _0}\left| \hat
u_n\left(t\right)\right|^2\lambda _0\left(\vartheta _0,t\right){\rm d}t<\delta  .
\end{align*}
\item[D]. We check \eqref{16} by direct calculations.

\end{description}

Having {\bf A}-{\bf C} by Theorem A.22 in \cite{IH81}  we obtain \eqref{15}.

To prove {\bf A}  we recall that by the central limit theorem
\begin{align}
\label{w}
W_n\left(t\right)=\sqrt{\frac{n}{\beta _0}}\left(\hat\Lambda
_n\left(t\right)-\beta _0\Lambda _0\left(\frac{t-\alpha_0 }{\beta
  _0}\right)\right) \Longrightarrow W\left(\Lambda _0\left(\frac{t-\alpha
  _0}{\beta _0}\right)\right) ,
\end{align}
where $W\left(r \right),0\leq r\leq \Lambda _0\left(\infty \right)$ is a
Wiener process.  Moreover, the vector $W_ {k,n} \equiv \left(W_n\left(t_1\right),\ldots,
W_n\left(t_k\right)\right)$ for any $k\geq 1$ and $t_i\in \R$ is
asymptotically normal
\begin{align*}
 W_ {k,n} \Longrightarrow W_ {k } =
  \left(W\left(\Lambda _0\left(\frac{t_1-\alpha
  _0}{\beta _0}\right)\right)  ,\ldots, W\left(\Lambda _0\left(\frac{t_k-\alpha
  _0}{\beta _0}\right)\right)   \right).
\end{align*}
 We know as well that the MLE $\hat\vartheta _n$ is asymptotically normal.
 The Wiener process $W\left(\cdot \right)$ and the Gaussian vector $\zeta $
 are correlated. To clarify this dependence and to prove the joint asymptotic
 normality of the MLE and of this vector we recall how the asymptotic
 normality of the MLE can be proved. We follow below the approach developed by
 Ibragimov and Khasminskii \cite{IH81}.

Introduce the normalized likelihood ratio
$Z_n\left(v\right)=\frac{L\left(\vartheta
  _0+\frac{v}{\sqrt{n}},X^n\right)}{L\left(\vartheta _0,X^n\right)} ,\;
  v\in \VV_n.$
Here $\VV_n=\left\{v\; :\;\vartheta _0+\frac{v}{\sqrt{n}}\in\Theta
\right\}$. Under the presented here conditions ${\scr R}$ the random field
$Z_n\left(v\right), v\in \VV_n $ admits the representation (LAN)
\begin{align}
\label{lan1}
\ln Z_n\left(v\right)=\langle v, S_n\left(\vartheta
_0,X^n\right)\rangle-\frac{1}{2} v^\tau \II\left(\vartheta _0\right)v+m_n,
\end{align}
where $m_n\rightarrow 0$ and the vector
\begin{align*}
S_n\left(\vartheta
_0,X^n\right)=\frac{1}{\sqrt{n}}\sum_{j=1}^{n}\int_{\R}^{}\frac{\dot \lambda 
  _0\left(\vartheta _0,t\right)}{\lambda
  _0\left(\vartheta _0,t\right)} \left[{\rm d}X_j\left(t\right)-\lambda
  _0\left(\vartheta_0 ,t\right){\rm d}t\right] 
\end{align*}
By the central limit theorem
\begin{align}
\label{lan2}
S_n\left(\vartheta _0,X^n\right)\Longrightarrow  S\left(\vartheta
_0\right)\sim {\cal N}\left(0,\II\left(\vartheta _0\right)\right).
\end{align}
Let us denote the limit random field
\begin{align*}
Z\left(v\right)=\exp\left\{\langle v, S\left(\vartheta
_0\right)\rangle-\frac{1}{2} v^\tau \II\left(\vartheta _0\right)v
\right\},\qquad v\in {\cal R}^2.
\end{align*}
Recall that we have the representation
\begin{align*}
 S\left(\vartheta _0\right)&=\sqrt{\beta _0}\int_{\R}^{}\frac{\dot \lambda
   _0\left(\frac{t-\alpha _0}{\beta _0}\right)}{\lambda _0\left(\frac{t-\alpha
     _0}{\beta _0}\right)}\;{\rm d}W\left(\Lambda _0\left(\frac{t-\alpha
   _0}{\beta _0}\right)\right)\\ 
&=\sqrt{\beta _0}\int_{\R}^{}\frac{\dot
   \lambda _0\left(s\right)}{\lambda _0\left(s\right)}\;{\rm d}W\left(\Lambda
 _0\left(s\right)\right)
\end{align*}
with the same Wiener process as in \eqref{w}.
Moreover, for the MLE we have the limit
\begin{align*}
\sqrt{\frac{n}{\beta _0}}\left(\hat\vartheta _n-\vartheta
_0\right)\Longrightarrow \zeta &=\II_*^{-1} \int_{\R}^{}\frac{\ell\left(s\right)}{\lambda
  _0\left(s\right)}\;{\rm d}W\left(\Lambda _0\left(s\right)\right),
\end{align*}
where the vector   $\ell\left(s\right)=-{\lambda
  _0'\left(s\right)}\left(1, s\right)^\tau $ (see \eqref{def}).
This representation, which we prove below, allows us to say what is the
correlation between $W\left(\Lambda _0\left(s\right)\right)$ and $\zeta $:
\begin{align*}
\Ex W\left(\Lambda _0\left(t\right)\right)\zeta &=\Ex \left[ W\left(\Lambda
  _0\left(t\right)\right)\II_*^{-1} \int_{\R}^{}\frac{ \ell
    _0\left(s\right)}{\lambda _0\left(s\right)}\;{\rm d}W\left(\Lambda
  _0\left(s\right)\right)\right]\\
&=\II_*^{-1} \int_{-\infty }^{t}\frac{ \ell
    _0\left(s\right)}{\lambda _0\left(s\right)}\;{\rm d}\Lambda
  _0\left(s\right)=\II_*^{-1} \int_{-\infty }^{t}{\ell
    _0\left(s\right)}\;{\rm d}s.
\end{align*}

Let us return to the proof of the asymptotic normality of the MLE. The random
field $Z_n\left(v\right),v\in\VV_n$ we extend on the whole plane $\R^2$
continuously decreasing to zero outside of $\VV_n$. Denote $\left({\cal
  C}_0\left(\R^2\right), {\scr B}\right)$ the measurable space of the
continuous 
random surfaces tending to zero at infinity  with the uniform metrics and Borelian $\sigma
$-algebra. Introduce the measures ${\bf Q}_n$ and ${\bf Q}$ induced by the
realizations of $Z_n\left(\cdot \right)$ and $Z\left(\cdot \right)$ in the
space $\left({\cal C}_0\left(\R^2\right), {\scr B}\right)$
respectively. Suppose that we already proved the weak convergence 
\begin{align}
\label{wc}
{\bf Q}_n\;\Longrightarrow \; {\bf Q}.
\end{align}
Then we have the convergence of the distributions of the continuous
functionals $\Psi \left(Z_n\right)$ to the distribution of $\Psi
\left(Z\right)$. 
 Consider a
convex set $\BB\in \R^2$. We can write 
\begin{align*}
&{\bf Q}_n\left(\sqrt{n}\left(\hat\vartheta _n-\vartheta _0\right)\in
\BB\right)\\
&\qquad\qquad ={\bf Q}_n \left(\sup_{\sqrt{n}\left(\vartheta-\vartheta _0\right)\in
\BB}L\left(\vartheta ,X^n\right)> \sup_{\sqrt{n}\left(\vartheta -\vartheta _0\right)\not\in
\BB}L\left(\vartheta ,X^n\right)\right)\\
&\qquad\qquad  ={\bf Q}_n \left(\sup_{\sqrt{n}\left(\vartheta-\vartheta _0\right)\in
\BB}\frac{L\left(\vartheta ,X^n\right)}{L\left(\vartheta_0 ,X^n\right) }>
\sup_{\sqrt{n}\left(\vartheta -\vartheta _0\right)\not\in 
\BB}\frac{L\left(\vartheta ,X^n\right)}{L\left(\vartheta_0
  ,X^n\right)}\right)\\
&\qquad \qquad =
{\bf Q}_n \left(\sup_{v\in
\BB}Z_n\left(v\right)> \sup_{v\not\in
\BB}Z_n\left(v\right)\right)\longrightarrow {\bf Q} \left(\sup_{v\in
\BB}Z\left(v\right)> \sup_{v\not\in
\BB}Z\left(v\right)\right)\\
&\qquad \qquad ={\bf Q}\left(\II\left(\vartheta _0\right)^{-1}S\left(\vartheta
_0\right)\in \BB\right) .
\end{align*}
Note that $\psi\left(z\right)=\sup_{v\in \BB}z\left(v\right)- \sup_{v\not\in
  \BB}z\left(v\right) $ is a continuous functional on the space $\left({\cal C}_0\left(\R^2\right), {\scr B}\right)$. The random function $Z\left(\cdot \right)$ takes its
maximum at the point $\hat v=\II\left(\vartheta _0\right)^{-1}S\left(\vartheta
_0\right) $. To prove the joint convergence in distribution of the vector
$W_{k,n}$ and $\hat v_n=\sqrt{\frac{n}{\beta _0}}\left(\hat\vartheta _n-\vartheta _0\right)$ we
denote $ R_n=\left(W_{k,n}, Z_n\left(v\right),v\in \R^2\right)$
introduce the product space ${\cal X}=\R^k\times {\cal C}_0\left(\R^2\right)$
with the corresponding Borelian $\sigma $-algebra ${\scr B}_*$. To verify the
weak convergence $R_n\Rightarrow R$, where $
R=\left(W_k,Z\left(v\right),v\in\R^2\right)$ we \\
a) prove the convergence of the
finite-dimensional distributions
\begin{align*}
\Bigl( W_{k,n}, Z_n\left(v_1\right),\ldots,
Z_n\left(v_m\right)\Bigr)\Longrightarrow  \Bigl( W_{k}, Z\left(v_1\right),\ldots,
Z_n\left(v_m\right)\Bigr)
\end{align*}
b) prove the tightness of the corresponding family of measures. 

The convergence a) follows from the LAN \eqref{lan1}, \eqref{lan2}. The prove
of b) is a part of the Theorem 1.10.1 in \cite{IH81}. The conditions ${\scr
  R}$ are sufficient for the verification of the conditions {\bf B1}-{\bf B5}
 of the Theorem 1.10.1 in \cite{IH81}. Therefore we obtain the joint
 asymptotic normality of the vector 
$\Bigl( W_{k,n}, \hat v_n\Bigr)\Longrightarrow  \Bigl( W_{k},\zeta  \Bigr).$

Hence we obtain the convergence of the finite-dimensional distributions \eqref{fdd}.
Let us check {\bf B}. We have
\begin{align*}
\hat u_n\left(t_1\right)-\hat
u_n\left(t_2\right)=W_n\left(t_1\right)-W_n\left(t_2\right)+\langle
\hat v_n ,\dot\Lambda (\vartheta_0 ,t_1)-\dot\Lambda (\vartheta_0
,t_2)\rangle. 
\end{align*}
Hence ($t_1<t_2$)
\begin{align*}
&\Ex_{\vartheta
  _0}\left|W_n\left(t_1\right)-W_n\left(t_2\right)\right|^2\\
&\qquad \quad =\Ex_{\vartheta
  _0}\left(\frac{1}{\sqrt{n}}
\sum_{j=1}^{n} \left[X_j\left(t_1\right)-X_j\left(t_2\right)-\Lambda
  _0\left(\vartheta _0,t_1\right)+\Lambda
  _0\left(\vartheta _0,t_2\right) \right]\right)^2\\
&\qquad \quad =\int_{t_1}^{t_2}\lambda _0\left(\vartheta _0,t\right)\,{\rm
  d}t\leq C\left|t_2-t_1\right|. 
\end{align*}
For the second term  we have
\begin{align*}
\Ex_{\vartheta  _0}\left|\langle
\hat v_n ,\dot\Lambda ( \vartheta_0 ,t_1)-\dot\Lambda ( \vartheta_0
,t_2)\rangle \right|^2&\leq \Ex_{\vartheta  _0}\left\|
\hat v_n \right\|^2\left\|
\dot\Lambda (\vartheta_0 ,t_1)-\dot\Lambda (\vartheta_0
,t_2) \right\|^2\\
&\leq C\,\left|t_2-t_1\right|^2\leq C\,L\left|t_2-t_1\right|.
\end{align*}

The inequality {\bf C} follows from the similar estimates.
\begin{align*}
&\int_{t>L}^{}\Ex_{\vartheta _0}\left|W_n\left(t\right)\right|^2\lambda
_0\left(\vartheta _0,t\right){\rm d}t =\int_{t>L}^{}\Lambda_0 \left(\vartheta
_0,t\right)\lambda _0\left(\vartheta _0,t\right){\rm d}t\\
&\qquad \qquad  <C\Lambda_0 \left(\infty \right) \left[\Lambda_0 \left(\infty \right)-\Lambda_0 \left(L \right) \right]\leq \delta 
\end{align*} 
because $\Lambda_0 \left(\infty \right)-\Lambda_0 \left(L \right)\rightarrow 0
$ as $L\rightarrow \infty $. 

The verification of {\bf D} easily follows from the given above estimates. 

Now the convergence $\hat\Delta _n\Rightarrow \Delta _0$ is proved as
follows. For any $d>0$ we take $L>0$ such that the estimate \eqref{dd}
holds.  The properties {\bf A}-{\bf C} according to the Theorem A.22 in
\cite{IH81} allow  us  to obtain the convergence 
\begin{align*}
\int_{-L}^{L}\hat u_n\left(t\right)^2\lambda _0\left(\vartheta
_0,t\right)\,{\rm d}t \Longrightarrow \int_{-L}^{L}\left[W\left(\Lambda
  _0\left(s\right)\right) -\langle \zeta
  ,V\left(s\right)\rangle\right]^2\lambda _0\left(s\right)\,{\rm d}s. 
\end{align*}
Therefore we proved that
\begin{align*}
\int_{\R}\hat u_n\left(t\right)^2\lambda _0\left(\vartheta
_0,t\right)\,{\rm d}t \Longrightarrow \int_{\R}\left[W\left(\Lambda
  _0\left(s\right)\right) -\langle \zeta
  ,V\left(s\right)\rangle\right]^2\lambda _0\left(s\right)\,{\rm d}s.
\end{align*}

Hence the test $\Psi_n\in {\scr K}_\varepsilon . $
\bigskip

There are several related problems which can be solved using such
approach. For example, in the case of periodic Poisson process with unknown
phase and frequency we have once more the model with shift and scale
parameters, but there is an essential difference too. The rate of convergence
of the estimate of the frequency is $n^{3/2}$.

 {\bf Acknowledgment.}  This work
was done under partial financial support of the grant of RSF 14-49-00079.

\end{document}